\documentclass[a4paper,12pt,reqno]{amsart}

\usepackage{amsmath,amssymb,amsthm,bm}
\usepackage{caption,graphicx}
\usepackage[english]{babel}
\usepackage{amscd}
\usepackage{amsgen}
\usepackage[final]{epsfig}
\usepackage{latexsym}
\usepackage{amsfonts}
\usepackage{url}
\usepackage{slashed}
\usepackage[all]{xy}
\usepackage{here}
\usepackage{comment}
\usepackage{fixmath}

\newtheorem{thm}{Theorem}
\newtheorem{lemm}{Lemma}

\newtheorem{defn}{Definition}

\allowdisplaybreaks

\begin{document}

\title{Conformal invariants of $3$-Braids and Counting Functions }

\author{Burglind J\"oricke }

\address{IHES, 35 Route de Chartres\\ 91440 Bures-sur-Yvette\\ France}

\email{joericke@googlemail.com}




\begin{abstract} We consider a conformal invariant of braids, the extremal length with totally real horizontal boundary values $\Lambda_{tr}$. The invariant descends to an invariant of elements of $\mathcal{B}_n\diagup\mathcal{Z}_n$, the braid group modulo its center. We prove that the number of elements of $\mathcal{B}_3\diagup\mathcal{Z}_3$ of positive $\Lambda_{tr}$ grows exponentially. The estimate applies to obtain effective finiteness theorems in the spirit of the geometric Shafarevich conjecture over Riemann surfaces of second kind. As a corollary we obtain another proof of the exponential growth
of the number of conjugacy classes of $\mathcal{B}_3\diagup\mathcal{Z}_3$ with positive entropy not exceeding $Y$.
\end{abstract}

\maketitle

\centerline \today

\medskip

In the paper \cite{Jo2} a conformal braid invariant is defined, the extremal length with totally real horizontal boundary values $\Lambda_{tr}$. Unlike the entropy the invariant $\Lambda_{tr}$ distinguishes in many cases the elements of a conjugacy class of braids.
Both invariants, $\Lambda_{tr}$ and the entropy, do not change under multiplication by an element of the center of the braid group, hence, they descend to invariants of $\mathcal{B}_n\diagup\mathcal{Z}_n$, the Artin braid group $\mathcal{B}_n$ modulo its center $\mathcal{Z}_n$.
We prove in this paper that the number of elements of $\mathcal{B}_3\diagup\mathcal{Z}_3$ with positive $\Lambda_{tr}$
not exceeding a positive number $Y$ grows exponentially.
As a corollary we obtain another proof of the exponential growth of the number of conjugacy classes of elements of $\mathcal{B}_3\diagup\mathcal{Z}_3$ that have positive entropy not exceeding  $Y$. Our proof does not use deep techniques from Teichm\"uller theory.

The first result that states exponential growth of the entropy counting function is due to Veech \cite{Ve}. More precisely, he
considered
conjugacy classes of pseudo-Anosov elements of the mapping class group of a closed  Riemann surface $S$, maybe with distinguished points $E_n$, with hyperbolic universal covering of $S\setminus E_n$, and proved that the number of classes  with entropy not exceeding a positive number $Y$, grows exponentially in $Y$.   Notice that the $3$-braids with positive entropy are exactly the $3$-braids that correspond to pseudo-Anosov elements of the mapping class group of the three-punctured complex plane.
The precise asymptotic for the entropy counting function is given by Eskin and Mirzakhani \cite{EsMi} for closed Riemann surfaces of genus at least $2$. Both papers, \cite{Ve} and \cite{EsMi}, use deep techniques of Teichm\"uller theory, in particular, they are based on the study of the Teichm\"uller flow.

The growth estimates for the counting function related to the extremal length allow to obtain effective finiteness theorems in the spirit of the Geometric Shafarevich Conjecture for the case when the base manifold is of second kind (\cite{Jo4}).

Before defining the extremal length with totally real horizontal boundary values we recall some facts.
We consider an $n$-braid as a homotopy class of loops with base point in the symmetrized configuration space. More detailed,
let $C_n(\mathbb{C}) \stackrel {def} = \{(z_1,\ldots,z_n) \in \mathbb{C}^n: z_j \neq z_k\, \mbox{if}\, j\neq k \}\,$ be the configuration space of $n$ particles moving in the complex plane without collision. The symmetric group $\mathcal{S}_n$   acts on  $C_n(\mathbb{C})$ .
The quotient $ C_n(\mathbb{C})\diagup \mathcal{S}_n$ is called the symmetrized configuration space. Its elements $E_n$ are unordered tuples of $n$ pairwise disjoint elements of $\mathbb{C}$, which are denoted by $E_n=\{z_1,\ldots,z_n\}$. We also consider $E_n$ as subset of
$\mathbb{C}$ consisting of $n$ points.

 We identify $\mathcal{B}_n$ with the fundamental group $\pi_1(C_n(\mathbb{C})\diagup \mathcal{S}_n,E_n)$ of the symmetrized configuration space with base point $E_n \in C_n(\mathbb{C})\diagup \mathcal{S}_n$.
Recall that there is an isomorphism  $\mathcal{B}_n \ni b \to \mathfrak{m}(b)\in \mathfrak{M}(\mathbb{D};\partial  \mathbb{D}, E_n)$  from the Artin braid group $\mathcal{B}_n$ of $n$-braids with base point $E_n$ onto the group $ \mathfrak{M}(\mathbb{D};\partial  \mathbb{D}, E_n)$ of isotopy classes of self-homeomorphisms of the closed disc $\overline{\mathbb{D}}$, which fix the boundary circle pointwise and the set $E_n$ setwise. The points of $E_n$ may be permuted. The pure braid group $\mathcal{PB}_n$ with base point $E_n$ is the group of braids that correspond to homeomorphisms which fix $E_n$ pointwise.
Further, there is a bijective correspondence between elements $\mathbold{b} \in
\mathcal{B}_n \diagup \mathcal{Z}_n$ and isotopy classes
$\mathfrak{m}_{\infty}(\mathbold{b})$ of orientation preserving self-homeomorphisms of the Riemann
sphere $\mathbb{P}^1$ with set of distinguished points $E_n \cup \{\infty\}$, more
precisely the homeomorphisms contained in $\mathfrak{m}_{\infty}(\mathbold{b})$ fix $\infty$
and map the set $E_n$ (the base point) onto itself (maybe, with a permutation of the
elements of $E_n$).

Notice that each self-homeomorphism of a punctured Riemann surface $S\setminus E_n$ extends to a self-homeomorphism of the closed surface $S$ that fixes the set of distinguished points $E_n$. We will identify self-homeomorphisms of punctured surfaces with self-homeomorphisms of closed surfaces with distinguished points.

The real configuration space $C_n (\mathbb {R})$ is defined in the same way as the complex configuration space
$C_n (\mathbb {C})$ with $\mathbb{C}$ replaced by $\mathbb{R}$.
The totally real subspace $C_n (\mathbb {R})
\diagup{\mathcal S}_n\,$ of the symmetrized configuration space  $\,C_n (\mathbb {C}) \diagup {\mathcal
S}_n\,$ is connected and simply connected.
Take a base point $E_n\in C_n
(\mathbb
{R}) \diagup {\mathcal S}_n $ . The Artin braid group $\mathcal{B}_n \cong \pi_1(\,C_n (\mathbb {C}) \diagup
{\mathcal S}_n\,,\; E_n\,)$
is isomorphic to the relative fundamental
group
$\pi_1(\,C_n(\mathbb {C}) \diagup {\mathcal S}_n\,,\; C_n
(\mathbb
{R}) \diagup {\mathcal S}_n \,)$. The elements of the latter
group
are homotopy classes of arcs in $\,C_n (\mathbb {C}) \diagup
{\mathcal S}_n\,$ with endpoints in $\,C_n(\mathbb {R})
\diagup
{\mathcal S}_n\,$.

The isomorphism between $\pi_1(\,C_n(\mathbb {C}) \diagup {\mathcal S}_n\,,\; C_n
(\mathbb{R}) \diagup {\mathcal S}_n \,)$ and the Artin braid group with generators $\sigma_j, \, j=1,\ldots,n-1,$ is given as follows. Take any arc $\gamma:[0,1]\to C_n(\mathbb {C}) \diagup {\mathcal S}_n$ with endpoints in $C_n(\mathbb
{R}) \diagup {\mathcal S}_n$ that represents an element of the relative fundamental group $\pi_1(\,C_n(\mathbb {C}) \diagup {\mathcal S}_n\,,\; C_n
(\mathbb{R}) \diagup {\mathcal S}_n \,)$. Associate to the arc $\gamma$ its lift  $\tilde {\gamma}=(\tilde {\gamma}_1,\ldots,\tilde {\gamma}_n): [0,1]\to C_n(\mathbb {C})$ for which $\tilde {\gamma}_1(0)< \ldots <\tilde {\gamma}_n(0)$. For $j=1,\ldots,n-1,$ the isomorphism maps the generator $\sigma_j$ to  the homotopy class of  arcs in  $C_n(\mathbb {C})\diagup \mathcal{S}_n$ with endpoints in $C_n(\mathbb{R})\diagup \mathcal{S}_n$, whose associated lift is represented by the following mapping. The mapping moves the points $\tilde {\gamma}_j(0)$ and $\tilde {\gamma}_{j+1}(0)$ with uniform speed in positive direction along a half-circle with center at the point $\frac{1}{2}(\tilde {\gamma}_j(0)+\tilde {\gamma}_{j+1}(0))$, and
fixes $\tilde {\gamma}_{j'}(0)$ for all $j'\neq j,j+1$.

Let $R$ be an open rectangle in the complex plane $\mathbb{C}$.
Speaking about rectangles we will always have in mind rectangles with sides
parallel to the coordinate axes. Denote the length of the horizontal
sides of $R$
by $\sf b$ and the length of the vertical sides by
$\sf a$.
(For instance, we may consider $R= \{z=x+iy:0<x<{\sf{b}},\,
0<y<{\sf{a}}\,\}$.)
The extremal length of $R$ introduced by Ahlfors \cite{A1}
equals $\lambda(R)=\frac{\sf a}{\sf
b}$.

Let $b \in \mathcal{B}_n$ be a braid.
Denote its image under the isomorphism from $\mathcal{B}_n$ to the relative fundamental group
$\pi_1(\,C_n
(\mathbb {C}) \diagup {\mathcal S}_n\,,\; C_n (\mathbb {R})
\diagup
{\mathcal S}_n\, )$ by $b_{tr}$. For a rectangle $R$ as above let $f:R
\to \,C_n
(\mathbb {C}) \diagup {\mathcal S}_n\,$ be a mapping which
admits a continuous extension to the closure $\bar R$ (denoted again by $f$) which
maps the (open) horizontal sides into $\,C_n (\mathbb {R})
\diagup
{\mathcal S}_n\,$. We say that the mapping represents $b_{tr}$
if for each maximal vertical line segment contained in $R$
(i.e. for $R$
intersected with any vertical line in $\mathbb{C}$) the
restriction of
$f$ to the closure of the line segment represents $b_{tr}$.

The definition of the  extremal length of an $n$-braid with
totally real horizontal boundary values is the following (\cite{Jo2}, \cite{Jo3}).

\begin{defn}\label{def1} Let $b \in \mathcal{B}_n$ be an
$n$-braid. The extremal length $\Lambda_{tr}(b)$ with totally
real horizontal
boundary values is defined as
\begin{align}
\Lambda_{tr}(b)=& \inf \{\lambda(R): R\, \mbox{ a rectangle
which
admits a holomorphic map to} \nonumber \\
&C_n (\mathbb {C}) \diagup {\mathcal S}_n \,\mbox{ that
represents}\; b_{tr}\}\,.\nonumber
\end{align}
\end{defn}
We recall that the centralizer $\mathcal{Z}_n$ of the Artin braid group $\mathcal{B}_n$  is the subgroup generated by the element $\Delta_n^2$ with $\Delta_n \stackrel{def}= (\sigma_1 \,\sigma_2 \ldots \sigma_{n-1})\,(\sigma_1 \,\sigma_2 \ldots \sigma_{n-2}) \, \ldots (\sigma_1 \,\sigma_2)\, \sigma_1$.
The invariant $\Lambda_{tr}$ does not change under multiplication by $\Delta_n$
(see e.g. Lemma 1 in \cite{Jo2}), and, hence, it is an invariant on $\mathcal{B}_n \diagup \mathcal{Z}_n$.

Our results concern the case $n=3$.  The bounds in the main theorems can be improved.
For the sake of simplicity of the proof we restrict ourselves to these estimates.
We consider first the
pure braid group modulo its center, ${\mathcal{PB}_3}\diagup \mathcal{Z}_3$, which is a free group in two generators $\sigma_j^2\diagup \mathcal{Z}_3$, $j=1,2$.
The counting function $N^{\Lambda}_{\mathcal{PB}_3}(Y),\, Y \in (0,\infty),$ is defined as follows. For each positive parameter $Y$ the value of $N^{\Lambda}_{\mathcal{PB}_3}(Y)$ is equal to the number of elements
$\mathbold{b} \in \mathcal{PB}_3\diagup \mathcal{Z}_3$
with $0 < \Lambda_{tr}(\mathbold{b})\leq Y$. Note that here we count elements of $\mathcal{PB}_3\diagup \mathcal{Z}_3$ rather than conjugacy classes of such elements. We wish to point out that (with the mentioned
choice of the generators) the condition $\Lambda_{tr}(\mathbold{b})>0$ excludes exactly the classes in  $\mathcal{PB}_3\diagup \mathcal{Z}_3$ of the even powers $\sigma_j^{2k}$ of the standard generators of the braid group $\mathcal{B}_3$ (see Theorem 1 in \cite{Jo2}).

\begin{thm}\label{thm10.1} For all positive numbers $Y\geq 600 \log 8$ the inequality
\begin{equation}\label{eq10.1}
\frac{1}{2} \exp(\frac{Y}{900}) \leq N^{\Lambda}_{\mathcal{PB}_3}(Y)\leq  \frac{1}{2}e^{6 \pi Y}
\end{equation}
holds. The upper bound is true for all $Y>0$.
\end{thm}

Consider arbitrary $3$-braids. The counting function $N^{\Lambda}_{\mathcal{B}_3}(Y),\, Y \in (0,\infty),$
is defined as the number of elements  $\mathbold{b}\in \mathcal{B}_3 \diagup \mathcal{Z}_3$
with $0 < \Lambda_{tr}(\mathbold{b})\leq Y$. Note that the condition $\Lambda_{tr}(\mathbold{b})>0$ excludes exactly the elements $\mathbold{b}\in \mathcal{B}_3 \diagup \mathcal{Z}_3$ that are represented by $b=\sigma_j^k \Delta_3^{\ell}$
for $j=1$ or $2$, and $\ell = 0$ or $1$. (See Lemma \ref{lem10.1'} and Theorem \ref{thm10.3} below.)

\begin{thm}\label{thm10.2} For any positive number $Y \geq 600 \log 8$ the inequality
\begin{equation}\label{eq10.12}
\frac{1}{2} \exp(\frac{Y}{900})\leq N^{\Lambda}_{\mathcal{B}_3}(Y)\leq 4 e^{6 \pi Y}
\end{equation}
holds. The upper bound is true for all $Y > 0$.
\end{thm}

The entropy  $h(b)$ of braids $b \in \mathcal{B}_n$
is defined as the infimum over the topological entropy of all homeomorphisms in the mapping class $\mathfrak{m}(b)$
corresponding to $b$. For the definition of the topological entropy and elementary properties see \cite{AKM}.
The entropy $h(b)$ does not change under conjugation and does not change under multiplication by an element of the center of the braid group. Moreover, the entropy of a braid $b$ is equal to the infimum of entropies of the elements of the mapping class $\mathfrak{m}_{\infty}(b)$ that is
associated to the image of $b$ in $\mathcal{B}_3\diagup \mathcal{Z}_n$ (for a proof see \cite{Jo1}). 
We will speak about the entropy of a braid, or about the entropy of an element of $\mathcal{B}_n\diagup \mathcal{Z}_n$ or about the entropy of the conjugacy class of such elements.

The entropy of mapping classes of compact surfaces of genus at least $2$ was first treated in the Asterisque volume \cite{FLP} dedicated to Thurston's theory of surface homeomorphisms.
Thurston gave a classification of mapping classes on closed or punctured Riemann surfaces  of genus $g$ with $n\geq 0$ punctures.
Here $2g-2 +n>0$, so that the universal covering of the Riemann surface
is the upper half-plane.
According to Thurston a finite non-empty set of mutually disjoint closed Jordan curves $C=\{ C_1 ,\ldots , C_{\alpha}\}$ on a closed or punctured Riemann surface $S$ is called admissible if no $C_i$ is homotopic to a point in $S$, or to a puncture, or to a $C_j$ with $i \ne j$. Thurston calls an isotopy class of self-homeomorphisms of $S$ reducible if there exists an admissible system of curves  $C=\{ C_1 , \ldots , C_{\alpha}\}$ on $S$ so that some homeomorphism in the class (and, hence each homeomorphism in the class) maps $C$ to an isotopic system
of curves. In this case we say that $C$ reduces the class.
If there is no such system the class is called irreducible.
Thurston proved that an isotopy class of homeomorphisms is irreducible if and only if it (more precisely, its extension to the closed Riemann surface) either contains a periodic homeomorphism or it contains a so called pseudo-Anosov homeomorphism. Periodic self-homeomorphisms of closed Riemann surfaces (maybe, with distinguished points)  have zero entropy, pseudo-Anosov
homeomorphisms have positive entropy.

We will also call a braid reducible, irreducible, periodic, or pseudo-Anosov if its image in $\mathcal{B}_n\diagup \mathcal{Z}_n$
corresponds to a class containing a reducible, irreducible, periodic, or pseudo-Anosov homeomorphism, respectively.

Notice that the reducible elements of $\mathcal{B}_3\diagup \mathcal{Z}_3$  are exactly the conjugates of powers of $\sigma_1 \diagup \mathcal{Z}_3$. It is known that their entropy equals zero (for a proof see e.g. \cite{Jo2} Theorem 2). Since periodic elements of $\mathcal{B}_3\diagup \mathcal{Z}_3$ have zero entropy, the pseudo-Anosov $3$-braids are exactly the $3$-braids of positive entropy.

The entropy counting function $N^{entr}_{\mathcal{B}_n}(Y),\, Y>0,$ for $n$-braids is defined 
as the number of conjugacy classes of pseudo-Anosov elements of $\mathcal{B}_n \diagup \mathcal{Z}_n$ with entropy not exceeding $Y$. For $3$-braids the value $N^{entr}_{\mathcal{B}_3}(Y),\, Y>0$, is also the number of  conjugacy classes of elements of $\mathcal{B}_3 \diagup \mathcal{Z}_3$ with positive entropy not exceeding $Y$.
The following theorem is a corollary of Theorem \ref{thm10.2}.

\begin{thm}\label{thm10.0} 
For any number $Y \geq 600 \pi \log 8 $ the estimate
\begin{equation}\label{eq10.0}
\frac{1}{2}\exp(\frac{Y}{900 \pi})    \leq N^{entr}_{\mathcal{B}_3}(Y)\leq 4e^{12 Y}
\end{equation}
holds. The upper bound holds for all positive $Y$.
\end{thm}

We will first prove Theorem \ref{thm10.1}. Notice that for each element $\mathbold{b}$ of
the pure braid group modulo its center $\mathcal{PB}_3 \diagup \mathcal{Z}_3$
there is a unique element $b\in\mathcal{PB}_3$ that represents $\mathbold{b}$ and can be written as reduced word in $\sigma_1^2$ and $\sigma_2^2$. Indeed, this representative $b$ of $\mathbold{b}$ is determined by the property that the linking number between the first and the third strand equals zero. Multiplying an arbitrary element $b \in \mathbold{b}$
with a suitable power of $\Delta_3^2$ one obtaines the required element. Assigning to each element $\mathbold{b}\in\mathcal{PB}_3 \diagup \mathcal{Z}_3$ the mentioned element we obtain an isomorphism to the free group in two generators $\sigma_1^2$ and $\sigma_2^2$, or, equivalently to
the fundamental group $\pi_1(\mathbb{C} \setminus\{-1,1\},0)$ of the twice punctured complex plane with base point $0$ with generators $a_1$ and $a_2$ which correspond to $\sigma_1^2$ and  $\sigma_2^2$, respectively. Representatives of $a_1$ surround $-1$ counterclockwise and representatives of $a_2$ surround $1$ counterclockwise. If no confusion arises we will identify an element of $\mathcal{PB}_3 \diagup \mathcal{Z}_3\cong \pi_1(\mathbb{C} \setminus\{-1,1\},0)$ with the corresponding reduced word in the generators $a_1$ and $a_2$,
or, eqivalently, with the corresponding reduced word in $\sigma_1^2$ and $\sigma_2^2$. 
The proof of Theorem \ref{thm10.1} uses the syllable decomposition of words $w \in \pi_1(\mathbb{C}\setminus \{-1,1\},0)$ introduced in \cite{Jo2}, \cite{Jo3}. It is defined as follows.
Write a word $w \in \pi_1(\mathbb{C}\setminus \{-1,1\},0)$ in reduced form
$w= a_{j_1}^{k_1}\, a_{j_2}^{k_2} \, \ldots\,,$ and call the $a_{j_l}^{k_l}$ the terms of $w$. The syllables are of two kinds. First, each term $a_{j_i}^{k_i}$ with $|k_i|\geq 2$ is a syllable (we call it a syllable of first type). Secondly, any maximal sequence of consecutive terms  $a_{j_i}^{k_i}$  for which $|k_i|=1$ and all $k_i$ have the same sign is a syllable (we call it a syllable of second kind). This gives a uniquely defined decomposition into syllables. The degree or length of the syllable is the sum of absolute values of the exponents of terms appearing in the syllable. We make the convention that the number of syllables of the identity equals zero. 

For a non-trivial word $w \in \pi_1(\mathbb{C}\setminus \{-1,1\},0) \cong \mathcal{PB}_3\diagup \mathcal{Z}_3$    we put $\mathcal{L}_-(w)\stackrel{def}= \sum \log(3 d_k)$ and $  \mathcal{L}_+(w) \stackrel{def}=  \sum \log(4d_k)$, where each sum runs over the degrees $d_k$ of all syllables of $w$.
The main ingredient of the proof
is the following lemma.

\begin{lemm}\label{lem10.1a}
Let $N^{\mathcal{L}_{-}}_{\mathcal{PB}_3}$ be the function whose value at any $Y >0$ is the number of reduced words  $w \in \pi_1(\mathbb{C}\setminus \{-1,1\},0)$, $w \neq \mbox{Id}$, for
which $\mathcal{L}_{-}(w)\leq Y$. The following inequality
\begin{equation}\label{eq10.1a'}
N^{\mathcal{L}_-}_{\mathcal{PB}_3}(Y) \leq \frac{1}{2}e^{3Y}
\end{equation}
holds.
\end{lemm}
We need some preparation for the proof of Lemma \ref{lem10.1a}.
Consider all finite tuples $(d_1,\ldots,d_j)$,
where $j\geq 1$ is any natural number (depending on the tuple) and the $d_k \geq 1$ are natural numbers.
Before proving Lemma \ref{lem10.1a} we estimate for $Y>0$ the number $N^*(e^Y)$ of different ordered tuples $(d_1, d_2,\ldots,d_j)$ (with varying $j$) that may serve as the degrees of the syllables of words (counted from left to right) with $\sum_1^j \log(3 d_k) \leq Y$.
Put $X=e^Y$. Then $N^*(X)$ is the number of distinct tuples
with $\prod (3d_k)\leq X$.

Fix a natural number $j$. Denote by $N_j^*(X),\, X \geq 1,$  the number of tuples $(d_1,\ldots,d_j)$ for which $\prod_1^j (3d_k) \leq X$. The number is not zero if and only if $j \leq \frac{\log X}{\log 3}$. For $X \geq 1$ the equality
\begin{equation}\label{eq10.4}
N^*(X) = \sum_{j \in \mathbb{Z}: \,1\leq j\leq  \frac{\log X}{\log 3}} N_j^*(X)
\end{equation}
holds.

Notice first that
$N_1^*(X)= [\frac{ X}{ 3}]$, where $[x]$ denotes the largest integer not exceeding the positive number $x$. Indeed we are looking for the number of $d_1$'s for which $1 \leq d_1 \leq \frac{ X}{ 3}$.

The value of $N_j^*$ for $j \geq 2$ is estimated by the following lemma.
Notice that the lemma holds also for $j=1$ if in the inequality
\eqref{eq10.7} we define $ 0! \stackrel{def}= 1$.

\begin{lemm}\label{lem10.2}
Let $j\geq 2$. Then $N_j^*(X)=0$ for $X< 3^j$. For  $X\geq 3^j$
\begin{equation}\label{eq10.7}
N_j^*(X) \leq \frac{1}{(j-1)!} \frac{1}{3} (\frac{2}{3})^{j-1} X \Big(\log \big(\frac{1}{3} (\frac{2}{3})^{j-1} X\big)\Big)^{j-1}.
\end{equation}
\end{lemm}

\noindent {\bf Proof.}
The number $N_2^*(X)$ is the number of tuples $(d_1,d_2)$ for which $3d_1 \cdot 3d_2 \leq X$. Since $d_1 d_2 \geq 1$, $N_2^*(X)=0$ for $X<3^2$.  If
$X\geq 3^2$ the inequality $d_1 \leq \frac{X}{9}$ holds, and for given $d_1$ the number $d_2$ runs through all natural numbers with $1\leq d_2 \leq \frac{X}{9d_1}$. Hence, for $X\geq 3^2$

\begin{equation}\label{eq10.5}
N_2^*(X) \leq  \sum_{k \in \mathbb{Z}: \, 1 \leq k \leq  \frac{X}{3^2}} \frac{X}{3^2 k}.
\end{equation}
Put $a=\frac{X}{3^2}$ and $k'= \frac{k}{a}$. Since for positive numbers $k'$ and $\alpha$ with $k'>\alpha$ the inequality $\frac{1}{k'}\leq \frac{1}{\alpha} \int_{k'-\alpha}^{k'} \frac{dx}{x}$ holds we obtain
\begin{align}\label{eq10.6}
N_2^*(X) \leq   \sum_{k' \in \frac{1}{a}\mathbb{Z}:\, \frac{1}{a} \leq k' \leq  1 } \frac{1}{k'}
\leq  2a \int_{\frac{1}{2a}}^1 \frac{dx}{x} = 2a \log(2a) =
\frac{1}{3} \cdot \frac{2}{3}X  \cdot \log(\frac{1}{3}\cdot \frac{2}{3}X ).
\end{align}

For $j\geq 3$ we provide induction using the following fact. Let $0<x<1$. Then for any positive integer $j$ the value
$(\frac{(-\log x)^j}{x})'= -\frac{(-\log x)^{j-1}}{x^2}(j-\log x) $ is negative.
Hence, for $k'\in (0,1)$ and $0<\alpha < k'$
\begin{equation}\label{eq10.6a}
\frac{1}{k'}(-\log(k'))^j  \leq \frac{1}{\alpha} \int_{k'-\alpha}^{k'} \frac{(-\log x)^j}{x}dx \,.
\end{equation}
We saw that the lemma is true for $j=2$. Suppose it is true for $j$. Prove that then it holds for $j+1$.  The number $N_{j+1}^*(X)$ is the number of tuples $(d_1,\ldots,d_{j+1})$ for which $3d_1 \cdot \ldots \cdot 3d_{j+1} \leq X$. Hence, $N_{j+1}^*(X)=0$ if $X <3^{j+1}$. If $X \geq 3^{j+1}$, then
$d_1 \leq \frac{X}{3^{j+1}}$ and for given $d_1$ the tuples $(d_2,\ldots,d_{j+1})$ run through all tuples with $3 d_2 \cdot \ldots \cdot 3 d_{j+1} \leq \frac{X}{3d_1}$. Hence, for $X \geq 3^{j+1}$

\begin{align}\label{eq10.8}
N_{j+1}^*(X) =  & \sum_{k \in \mathbb{Z}:\, 1 \leq k \leq  \frac{X}{3^{j+1}}} N_j^*(\frac{X}{3k})\nonumber\\
\leq & \sum_{k \in \mathbb{Z}:\, 1 \leq k \leq  \frac{X}{3^{j+1}}} \frac{1}{(j-1)!}  \frac{1}{3} (\frac{2}{3})^{j-1}  \frac{X}{3k} \Big( \log\big(\frac{1}{3} (\frac{2}{3})^{j-1}  \frac{X}{3k}\big)\Big)^{j-1}.
\end{align}

Put $a= \frac{1}{9} (\frac{2}{3})^{j-1} X$ and $k'= \frac{k}{a}$. Then
\begin{align}\label{eq10.9}
N_{j+1}^*(X) \leq  & \frac{1}{(j-1)!}\sum_{k' \in \frac{1}{a}\mathbb{Z}:\frac{1}{a} \leq k' \leq  \frac{1}{2^{j-1}}} \frac{1}{k'} \big(\log(\frac{1}{k'})\big)^{j-1}\nonumber \\
\leq &  \frac{1}{(j-1)!} 2a \int_{\frac{1}{2a}}^{\frac{1}{2^{j-1}}} \frac{1}{x} (-\log x )^{j-1} dx
\leq   \frac{1}{(j-1)!} 2a \int_{\frac{1}{2a}}^1 \frac{1}{x} (-\log x )^{j-1} dx
\nonumber\\
= & (-1)^{j-1} \frac{1}{(j-1)!} 2a \frac{1}{j}(\log x)^j\mid_{\frac{1}{2a}}^{1}.
\end{align}
We obtain
\begin{align}\label{eq10.10}
N_{j+1}^*(X)
\leq & \frac{1}{j!} 2a \big(\log (2a)\big)^j = \frac{1}{j!}\frac{1}{3} (\frac{2}{3})^{j}X \Big(\log\big( \frac{1}{3} (\frac{2}{3})^{j}X\big)\Big)^j.
\end{align}
Lemma \ref{lem10.2} is proved.
\hfill $\Box$

\bigskip

Lemma \ref{lem10.2} implies the following upper bound for $N^*$.

\begin{lemm}\label{lem10.1}
For $X<3$ the function $N^*(X)$ vanishes. Moreover, for any positive number X
\begin{equation}\label{eq10.0}
N^*(X) \leq (\frac{X}{3})^{\frac{5}{3}}.
\end{equation}
\end{lemm}

\noindent {\bf Proof of Lemma \ref{lem10.1}.}
Since all $N_j^*(X)$ vanish for $X<3$ the value $N^*(X)$ vanishes for such $X$.
For $X\geq 3$ equation \eqref{eq10.4} and Lemma \ref{lem10.2} imply
\begin{align}\label{eq10.11}
N^*(X) \leq & \sum_{j \in \mathbb{Z}; \, 1 \leq j \leq \frac{\log X}{\log 3}}  \frac{1}{(j-1)!} \frac{1}{3} (\frac{2}{3})^{j-1} X \Big(\log\big(\frac{1}{3} (\frac{2}{3})^{j-1} X\big)\Big)^{j-1} \nonumber \\
\leq & \; \frac{1}{3} X \sum_{j \in \mathbb{Z}: \, 1 \leq j \leq \frac{\log X}{\log 3}} \frac{1}{(j-1)!} \big(\frac{2}{3} \log (\frac{1}{3}X)\big)^{j-1} \nonumber \\  \leq & \; \frac{1}{3} X \exp\big(\frac{2}{3}\log(\frac{1}{3}X)\big)=
(\frac{X}{3})^{\frac{5}{3}}.
\end{align}
Lemma \ref{lem10.1} is proved. \hfill $\Box$

\bigskip

\noindent {\bf Proof of Lemma \ref{lem10.1a}.}
We assume that $[\frac{Y}{\log 3}]\geq 1$. Otherwise $N^*(e^Y)$ vanishes, and therefore $N^{\mathcal{L}_-}_{\mathcal{PB}_3}(Y)=0$, and the inequality is satisfied.
We will use the notation $N_j^{\mathcal{L}_-}(Y), \, j\geq 1,$ for the number of different reduced words $w$ in $\pi_1(\mathbb{C} \setminus \{-1,1\}, 0)$ that  consist of $j$ syllables and satisfy the inequality $\mathcal{L}_-(w)= \sum_{k=1}^j \log(3d_k) \leq  Y$.
Then $N^{\mathcal{L}_-}_{\mathcal{PB}_3}(Y)=  \sum_{j=1}^{j_0}N_j^{\mathcal{L}_-}(Y)$ with $j_0\stackrel{def}=[\frac{Y}{\log 3}]$.
We will estimate $N_j^{\mathcal{L}_-}(Y)$ by $N_j^*(X)$ with $X= e^{ Y}$. Recall that $N^*(e^Y)$ is the number of different tuples $(d_1,\ldots,d_j)$ with $d_k \geq 1$ for which $\prod_{k=1}^j (3d_k)\leq e^{ Y}$.

For this purpose we take a tuple $(d_1,\ldots,d_j)$ and estimate the number of different reduced words with tuple of lengths of syllables (from left to right) equal to  $(d_1,\ldots,d_j)$. The first syllable can be of type (1) or of type (2). The type of the syllable and the first letter of the syllable (which may be $a_1^{\pm 1}$ or $a_2^{\pm 1}$) together with its length completely determine the syllable. Hence, there are at most $8$ different choices for the first syllable if we require the syllable to have exactly degree $d_1$.
For all other syllables
the first letter of the syllable cannot be $a_i^{\pm 1}$ if the last letter in the preceding syllable is $a_i^{\pm 1}$ for the same $a_i$. Hence, for all but the first syllable there are at most $4$ choices given the degree of the syllable and the preceding syllable.

We showed that for all $j=1,\ldots, j_0 =[\frac{Y}{\log 3}]$, and each tuple $(d_1,\ldots,d_j)$,  there are at most $2 \cdot 4^j$ different reduced words with tuple of lengths of syllables  equal to  $(d_1,\ldots,d_j)$.
Hence, for $Y  \geq \log 3$ the number $N^{\mathcal{L}_-}_{\mathcal{PB}_3}(Y)$ of reduced words
$w \in   \pi_1(\mathbb{C}\setminus \{-1,1\},0)$, $w\neq \mbox{Id}$,  with $\prod_1^{j_0} (3d_k) \leq \exp(Y)$ equals
\begin{equation}\label{eq10.3}
N^{\mathcal{L}_-}_{\mathcal{PB}_3}(Y)= \sum_{j=1}^{j_0}N_j^{\mathcal{L}_-}(Y)\;\,\leq\;\;
\sum_{j=1}^{j_0}
2\cdot 4^{j_0}N_j^*(e^{Y})\, =\,
2\cdot 4^{j_0} \cdot N^*(e^{Y}).
\end{equation}

Using the inequality $j_0 
\leq \frac{Y}{\log 3}$ and Lemma \ref{lem10.1} with $X=e^{Y}$ we obtain
the requested estimate by the value
\begin{equation}\label{eq10.3a}
2\cdot 3^{-\frac{5}{3}} \exp((\frac{ \log 4}{\log 3}+\frac{5}{3}) Y)< \frac{1}{2}\exp(3 Y).
\end{equation}
We used the inequalities $(\frac{ \log 4}{\log 3}+(\frac{5}{3}))<2.93$ and $4 \cdot 3^{-\frac{5}{3}} < 0.65<1$.
Lemma \ref{lem10.1a} is proved.
\hfill $\Box$

\bigskip

\noindent {\bf Proof of 
Theorem \ref{thm10.1}.}  We need the following inequality (see Theorem 1 in \cite{Jo2}) 
\begin{align}\label{eq10.2'}
\frac{1}{2 \pi} \mathcal{L}_-(w)
\leq \Lambda_{tr}(w) \leq 300 \mathcal{L}_+(w)
\end{align}
which holds for all reduced words $w$ representing elements in $\mathcal{PB}_3\diagup \mathcal{Z}_3\cong \pi_1(\mathbb{C}\setminus \{-1,1\},0)$ that are not equal to a power of $a_1$ or of $a_2$  or to the identity (equivalently, for which $\Lambda_{tr}(w)>0$).

The inequality \eqref{eq10.2'} implies the inclusion
$\{w:0<\Lambda_{tr}(w) \leq Y\} \subset \{w \neq \mbox{Id}: \mathcal{L}_-(w) \leq 2\pi Y \}$. We obtain
the inequality
$N_{\mathcal{PB}_3}^{\Lambda}(Y) \leq N^{\mathcal{L}_-}_{\mathcal{PB}_3}(2 \pi Y)$  and by
Lemma \ref{lem10.1a} the right hand side of this inequality does not exceed $\frac{1}{2}e^{6\pi Y}$. This gives the upper bound.

The lower bound is obtained as follows.
Consider all reduced words in $\pi_1(\mathbb{C}\setminus\{-1,1\},0)$ of the form
\begin{equation}\label{eq10.2''}
a_1^{2k_1} a_2^{2k_2} \ldots
\end{equation}
where each $k_i$ is equal to $1$ or $-1$. If $j$ is the number of syllables (i.e the number of  the $a_i^{k_i}$) of a word $w$ of the form \eqref{eq10.2''}, then
$\Lambda_{tr}(w) \leq 300 \mathcal{L}_+(w) =300  j \log 8$.
Consider the words of the mentioned form for which
$j=j_0 \stackrel{def}= [\frac{Y}{300 \log 8}]$. Since $j_0$ must be at least equal to $2$ we get the condition $Y \geq 600 \log 8$. For the chosen $j_0$  the extremal length of the considered words does not exceed $Y$. The number of different words of such kind is
$2^{j_0}= 2^ {[\frac{Y}{300 \log 8}]}\geq \exp(\log 2 \cdot( \frac{Y}{300 \log 8}-1))= \frac{1}{2} \exp( Y \frac{\log 2}{300\cdot 3 \log 2})= \frac{1}{2} \exp(\frac{Y}{900} )$ .
Theorem \ref{thm10.1} is proved. \hfill $\Box$

\medskip

Consider now arbitrary elements of the braid group modulo its center $\mathcal{B}_3\diagup \mathcal{Z}_3$ and their extremal length with totally real horizontal boundary values.

We need the following lemma and theorem
from \cite{Jo2} which are formulated in terms of braids.
\begin{lemm}\label{lem10.1'}{\rm{(\cite{Jo2}, Lemma 2)}}
Any braid $b\in \mathcal{B}_3$ which is not a power of $\Delta_3$ can be written in a unique way in the form
\begin{equation}\label{eq2'}
\sigma_j^k \, b_1 \, \Delta_3^{\ell}
\end{equation}
where $j=1$ or $j=2$, $k\neq 0$ is an integer, $\ell$ is a (not necessarily even) integer, and $b_1$ is a word in $\sigma_1^2$ and $\sigma_2^2$ in reduced form. If
$b_1$ is not the identity, then the first term of $b_1$ is a non-zero even
power of $\sigma_2$ if $j=1$, and the first term of $b_1$ is a non-zero even  power of
$\sigma_1$ if $j=2$.
\end{lemm}

For an integer $l\neq 0$ we put $q(l)=l$ if $l$ is even, and for odd $l$ we
denote by $q(l)$ the even integer neighbour of $l$ that is closest to zero. In other words, $q(l)=l$ if $l\neq 0$ is even, and for each odd integer $l\,,$  $q(l)= l
-{\mbox{sgn}}(l)$, where  ${\mbox{sgn}}(l)$
for a non-zero integer number $l$ equals $1$ if $l$ is positive,
and $-1$ if $l$ is negative. For a braid in form \eqref{eq2'} we put
$\vartheta(b) \stackrel{def}{=}\sigma_j^{q(k)} \, b_1$. $\vartheta(b)$ can be written as a word in $\sigma_1^2$ and $\sigma_2^2$. 
It is clear that $\vartheta(b \Delta_3)=\vartheta(b)$ for $b\in \mathcal{B}_3$.
For $\mathbold{b}\in \mathcal{B}_3\diagup \mathcal{Z}_3$ the value $\vartheta(\mathbold{b})$ is well-defined by the relation $\vartheta(\mathbold{b})\stackrel{def}= \vartheta(b)$ for any $b\in \mathcal{B}_3$ representing $\mathbold{b}$.
The following theorem holds.

\begin{thm}\label{thm10.3}{\rm{(\cite{Jo2}, Theorem 3)}}
Let $b \in \mathcal{B}_3$ be a (not necessarily pure) braid which is not a power of $\Delta_3$, and let $w$ be the reduced word 
in $\sigma_1^2$ and $\sigma_2^2$, that represents $\vartheta(b)$.
Then
$$
\frac{1}{2\pi}\mathcal{L}_-(w) \leq \Lambda_{tr}(b) \leq 300 \cdot
\mathcal{L}_+(w) \,,
$$
except in the case when $b=\sigma_j^{k}\,\Delta_3^{\ell}$, where $j=1$ or  $j=2$, $k\neq 0$ is an integral number, and $\ell$ is an arbitrary integer. In this case $\Lambda_{tr}(b)=0$.
\end{thm}

\medskip

\noindent {\bf Proof of Theorem \ref{thm10.2}.}
Take any element of $\mathcal{PB}_3 \diagup \mathcal{Z}_3$. Choose its unique representative that can be written as a reduced word $w$ in $\sigma_1^2$ and $\sigma_2^2$.
We describe now all elements $\mathbold{b}$
of $\mathcal{B}_3 \diagup \mathcal{Z}_3$ with
$\vartheta(\mathbold{b})=w$.
If $w\neq \mbox{Id}$ these are the elements represented by the following braids.
If the first term of $w$ is 
$\sigma_j^{2k}$ with $k\neq 0$, then the possibilities are ${b}=w \Delta_3^{\ell}$ with $\ell=0$ or $1$,  ${b}=\sigma_j^{\begin{tiny}{\mbox{sgn}}\end{tiny} k} w \Delta_3^{\ell}$ with $\ell=0$ or $1$, or  ${b}=\sigma_{j'}^{\pm 1} w \Delta_3^{\ell}$ with $\ell=0$ or $1$ and $\sigma_{j'}\neq
\sigma_{j}$. Hence, for $w\neq \mbox{Id}$ there are $8$ possible choices of elements $\mathbold{b}\in \mathcal{B}_3 \diagup \mathcal{Z}_3$ with $\vartheta(\mathbold{b})=w$.
By Theorem \ref{thm10.3} the set of $\mathbold{b} \in \mathcal{B}_3 \diagup \mathcal{Z}_3$ with $0 < \Lambda_{tr}(\mathbold{b}) \leq Y$ is contained in the set of  $\mathbold{b} \in \mathcal{B}_3 \diagup \mathcal{Z}_3$ with $\vartheta({\mathbold{b}}) = w \neq \mbox{Id}$, $ \mathcal{L}_-(w) \leq 2\pi Y$. We obtain
\begin{equation}\label{eq10.13}
N^{\Lambda}_{\mathcal{B}_3}(Y)\leq  8 N^{\mathcal{L}_-}_{{\mathcal{PB}}_3}(2\pi Y) .
\end{equation}
By Lemma \ref{lem10.1a} we obtain  $N^{\Lambda}_{\mathcal{B}_3}(Y)\leq 4 e^{6 \pi Y}$.

Since each pure $3$-braid is also an element of the braid group $\mathcal{B}_3$ the lower bound of Theorem \ref{thm10.1} provides also a lower bound for Theorem \ref{thm10.2}. Theorem \ref{thm10.2} is proved. \hfill $\Box$

\bigskip
We prepare the proof of Theorem \ref{thm10.0}.
A reduced word $w\neq \mbox{Id}$ representing an element of $\mathcal{PB}_3 \diagup\mathcal{Z}_3$ is called cyclically reduced, if either the word consists of a single term, or it has at least two terms and the first and the last term of the word are powers of different generators. Each reduced word which is not the identity is conjugate to a cyclically reduced word.

A cyclically reduced word $w$ representing an element of $\mathcal{PB}_3 \diagup\mathcal{Z}_3$ is called cyclically syllable reduced, if either the word consists of a single term, or all terms enter with equal power $+1$ or $-1$ (in these two cases the word consists of a single syllable), or the first and the last term of the word do not enter with equal power $+1$ or $-1$ (in this case the word contains at least two syllables). Each cyclically reduced word which is not the identity is conjugate to a cyclically syllable reduced word.

For the proof of theorem \ref{thm10.0} we will use the following theorem from \cite{Jo2}. 
\medskip

\begin{thm}\label{thm10.1a}{\rm{(\cite{Jo2}, Theorem 2)}}
Let $\mathbold{b}\in \mathcal{PB}_3\diagup \mathcal{Z}_3$ be represented by a cyclically syllable reduced word $w$
consisting of more than one syllable. Then
\begin{align}\label{eq10.2}
\frac{1}{2 \pi} \mathcal{L}_-(w)
\leq \frac{2}{\pi} h(\hat{\mathbold{b}}) \leq 300 \mathcal{L}_+(w).
\end{align}
\end{thm}

The extremal length in the sense of Ahlfors of a round annulus $A=\{z \in \mathbb{C}:  r_1 <|z|<r_2\},\;$ $0\leq r_1 < r_2 \leq \infty,\;$  equals $\lambda(A)=\frac{2\pi}{\log(\frac{r_2}{r_1})}$. There is a bijective correspondence between conjugacy classes $\hat b$ of $n$-braids and free isotopy classes of loops in $C_n(\mathbb{C})\diagup \mathcal{S}_n$.

A continuous mapping $f$ of an annulus $A= \{z \in \mathbb{C}:
\, r_1<|z|<r_2\},\;$ $0\leq r_1 < r_2 \leq \infty,\;$ into
$C_n(\mathbb{C}) \diagup \mathcal{S}_n$ represents a conjugacy class
$\hat b$ of $n$-braids if for
each circle $\,\{|z|=\rho \} \subset A\,$
the
loop $\,f:\{|z|=\rho \} \rightarrow  C_n(\mathbb{C}) \diagup
\mathcal{S}_n
\,$
represents the conjugacy class $\,\hat b$.
The extremal length $\Lambda(\hat
b)$ of $\hat b$ is defined as $ \Lambda(\hat b)= \inf_{A \in
\mathcal{A}}\,
\lambda(A),$ where
$\mathcal{A}$ denotes the set of all annuli which admit a
holomorphic mapping into
$C_n(\mathbb{C}) \diagup \mathcal{S}_n$ that represents $\hat b$.
By \cite{Jo} and \cite{Jo1} the entropy of a conjugacy class of $n$-braids $\hat b$ equals $h(\hat{b})= \frac{\pi}{2}\Lambda(\hat{b})$. The equality $h(\hat{\mathbold{b}})=h(\hat{b})$
holds for any braid $b$ representing $\mathbold{b}$ (see e.g. \cite{Jo}, \cite{Jo1}).
\medskip

\noindent {\bf Proof of the upper bound of Theorem \ref{thm10.0}.}
We represent each 
conjugacy class $\hat{\mathbold{ b}}$ of elements of $\mathcal{B}_3\diagup\mathcal{Z}_3$ by a conjugacy class $\hat b$ of elements of $\mathcal{B}_3$.
To each  conjugacy class $\hat b$ of elements of $\mathcal{B}_3$  with $h(\hat{b})>0$ and each positive number $\varepsilon$ we will associate a braid $b\in \mathcal{B}_3$ that represents  $\hat{ b}$
such that the inequality
\begin{equation}\label{eq10.30}
\Lambda_{tr}(b)< \frac{2}{\pi}   h(\hat{b})+\varepsilon
\end{equation}
holds.
For this purpose we represent the conjugacy class $\hat b$ of $b$ by a holomorphic map $g:A\to C_3(\mathbb{C})\diagup \mathcal{S}_3$ from an annulus $A$ of extremal length
\begin{equation}\label{eq10.29}
\lambda(A)<\frac{2}{\pi}h(\hat{b})+\varepsilon
\end{equation}
to the symmetrized configuration
space. By Lemma 6 of \cite{Jo2} each loop that represents $\hat b$ intersects the smooth real hypersurface
\begin{align}
\mathcal{H} \stackrel{def}= & \{ \{z_1,z_2,z_3\} \in C_3(\mathbb{C})\diagup \mathcal{S}_3: \mbox{the three points}\; z_1,z_2, z_3  \nonumber\\
& \mbox{are contained in a real line in the complex plane}\}
\end{align}
of $C_3(\mathbb{C})\diagup \mathcal{S}_3$.
By the Holomorphic Transversality Theorem \cite{KZ} we may assume, after shrinking $A$ (keeping inequality \eqref{eq10.29} ) and approximating $g$, that $g$ is holomorphic in a neighbourhood of the closure $\bar A$ of $A$ and is transversal to $\mathcal{H}$.
Hence, $L\stackrel{def}=\{z \in \bar{A}: g(z) \in \mathcal{H}\}$ is a smooth
real submanifold of $\bar{A}$ of real dimension $1$. Moreover, $L$ contains an arc $L_0 \subset A$ with endpoints on different boundary circles of $A$.

The set $A\setminus L_0$ is a curvilinear rectangle. By this we man that the set admits a conformal mapping  $\omega$ onto a rectangle $R=\{z=x+iy\in \mathbb{C}:\,0<x<1,\,0<y<\sf{a}\}$ of extremal length $\sf{a}$ with the following properties.
Consider a lift $(A\setminus L_0)\;\tilde{}\;$ of the set $A\setminus L_0$ to the universal covering of $\bar A$ and denote the projection map by $p$.
The lift of the mapping $\omega$ extends to a homeomorphism from the closure of
$(A\setminus L_0)\;\tilde{}\;$ onto the closure $\bar{R}$ of $R$ such that
the two components 
of $p^{-1}(\partial A \setminus \overline{L_0}) $  are mapped to the open vertical sides of the rectangle and the two components
of $p^{-1}(L_0)$ are mapped to the open horizontal sides of the rectangle.

The extremal length $\lambda(A\setminus L_0)$ is defined as the extremal length $\sf{a}$ of the rectangle $R$. The inequality $\lambda(A\setminus L_0)\leq \lambda(A) $ holds (see \cite{A1} , or inequality (5) of \cite{Jo4}, or the remark after the proof of Theorem
3 of \cite{Jo2}).


We may assume that $L_0$ contains a point $z_0$ for which $g(z_0)$ is contained in the real subspace $C_3(\mathbb{R})\diagup \mathcal{S}_3$, and, moreover, for some label $\{g_1(z_0),g_2(z_0),g_3(z_0)\}$ of the points of $g(z_0)$ the equalities $g_1(z_0)=0$, $g_3(z_0)=1$ hold.
To see this we
will use the following notation. Let $\mathfrak{A}$ be a complex affine mapping, i.e. $\mathfrak{A}(\zeta)=\alpha\,\zeta +\beta,\,\zeta \in \mathbb{C},$ where $\alpha $ and $\beta$ are complex numbers, $\alpha \neq 0$. For a point $E\in \mathcal{C}_3(\mathbb{C})\diagup \mathcal{S}_3$ we denote by
$\mathfrak{A}( E)$ the triple of points in $\mathbb{C}$ that is obtained by applying $\mathfrak{A}$ to each of the three points of $E$.

Recall that $g(L_0)\subset \mathcal{H}$. Take any point $z_0\in L_0 \cap A$.
Label the points of $g(z_0)$ in any way by $g_1(z_0),\,g_2(z_0),\,$
and $g_3(z_0)$.
We denote by $\mathfrak{A}_{z_0}$ the complex affine mapping $\zeta\to \frac{\zeta-g_1(z_0)}{g_3(z_0)-g_1(z_0)}, \, \zeta \in \mathbb{C}$. Consider the mapping $z \to \mathfrak{A}_{z_0}( g(z)),\, z \in A,$ which assigns to each point $z\in A$ the point in  $C_3(\mathbb{C})\diagup \mathcal{S}_3$ that is obtained by applying the complex affine mapping $\mathfrak{A}_{z_0}$ to $g(z)$. Then $\mathfrak{A}_{z_0}(g(z_0))=\{0, \frac{g_2(z_0)-g_1(z_0)}{g_3(z_0)-g_1(z_0)},1\}$  is contained in $\mathcal{C}_3(\mathbb{R})\diagup \mathcal{S}_3$, and has the required form.

Since the set of complex affine mappings (with the topology of $\mathbb{C}^*\times \mathbb{C}$) is connected, the mapping $\mathfrak{A}_{z_0} g$ is free isotopic to $g$ on $A$, hence also represents $\hat b$.
Replacing if necessary $g$ by $\mathfrak{A}_{z_0} g$, we assume now that $g(z_0) \in C_3(\mathbb{R})\diagup \mathcal{S}_3$, and has the required form.

Take the lift of the mapping $g\mid A \setminus L_0$ to a mapping 
$(g_1,g_2,g_3):A\setminus L_0 \to \mathcal{C}_3(\mathbb{C})$
for which the continuous extension $(g_1(z_0),\,g_2(z_0),g_3(z_0))$ to $z_0$ in clockwise direction satisfies $g_1(z_0)=0$ and $g_3(z_0)=1$.
We consider the mapping $z\to \alpha(z)\zeta+\beta(z) \stackrel{def}= \mathfrak{A}_{z}(g(z)),\, z \in A\setminus L_0,$ which is obtained by applying the variable complex affine mapping $\mathfrak{A}_z$,  $\mathfrak{A}_z(\zeta)= \frac{\zeta-g_1(z)}{g_3(z)-g_1(z)},$ to $g(z)$. Use the notation
$g_{\mathfrak{A}}(z)$ for the mapping $\mathfrak{A}_{z}(g(z)),\, z\in A\setminus L_0$. The value of this mapping for each $z\in A\setminus L_0$  is the triple $\{0,\frac{g_2(z)-g_1(z)}{g_3(z)-g_1(z)},1\}$. Hence, the continuous extension of this mapping to each of the strands of $L_0$
takes values in the real subspace
$C_3(\mathbb{R})\diagup \mathcal{S}_3$ of the symmetrized configuration space, and therefore the holomorphic mapping $g_{\mathfrak{A}}\circ \omega^{-1}$ on $R$ represents an element of the relative fundamental group $\pi_1(C_3(\mathbb{C})\diagup \mathcal{S}_3,C_3(\mathbb{R})\diagup \mathcal{S}_3)$.

Take a curve $\gamma:[0,1]\to A$ with $\gamma(0)=\gamma(1)=z_0$  that intersects $L_0$ only at $z_0$ and represents the counterclockwise generator of the fundamental group of $A$.
The mapping $g\circ\gamma$ represents a braid $b\in \hat b$ and also an element in the relative fundamental group $\pi_1(C_3(\mathbb{C})\diagup \mathcal{S}_3,C_3(\mathbb{R})\diagup \mathcal{S}_3)$.
The mapping $g_{\mathfrak{A}}\circ \gamma\mid (0,1)$ extends continuously to the closed interval $[0,1]$. The extension (denoted by $g_{\mathfrak{A}}\circ \gamma$)
represents an element of the relative fundamental group
$\pi_1(C_3(\mathbb{C})\diagup \mathcal{S}_3,C_3(\mathbb{R})\diagup \mathcal{S}_3)$. By the construction of $g_{\mathfrak{A}}$ the equality $g_{\mathfrak{A}}\circ \gamma(0)= g\circ \gamma (0)$ holds, and
the two elements of the relative fundamental group
represented by $g\circ \gamma$ and by $g_{\mathfrak{A}}\circ \gamma$ differ by a finite number of half-twists.
For the conformal mapping $\omega:A\setminus L_0\to R$ and the extremal length $\sf{a}$ of the rectangle $R$
the map $z \to e^{\frac{\pi}{\sf{a}}\omega(z)}, \, z\in A\setminus L_0,$ represents a half-twist. For $k\in \mathbb{Z}$ we consider the mapping
$ g_{k,\omega}(z) \stackrel{def}= e^{k\frac{\pi}{\sf{a}}\omega(z)} g_{\mathfrak{A}}(z)$ from $A\setminus L_0$ to $C_3(\mathbb{C})\diagup \mathcal{S}_3$. It has totally real horizontal boundary values.
There exists an integer number $k$ such that
$g\circ \gamma$ and $g_{k,\omega}\circ\gamma$
represent the same element of $\pi_1(C_3(\mathbb{C})\diagup \mathcal{S}_3,C_3(\mathbb{R})\diagup \mathcal{S}_3)$.
Hence, the mapping $g_{k,\omega}\circ \omega^{-1}:R\to C_3(\mathbb{C})\diagup\mathcal{S}_3$  represents the element $b_{tr}$ of the relative fundamental group that corresponds to the braid $b\in \hat{b}$ represented by $g\circ \gamma$. Since $\lambda(R)={\sf{a}}=\lambda(A\setminus L_0)\leq \lambda(A)$ we obtained
$\Lambda_{tr}(b)\leq \lambda(A)< \frac{2}{\pi}h(\hat{b})+\varepsilon$.

We achieved the following.
For each conjugacy class $\hat b$ of $3$-braids with $h(\hat{b})>0$ we obtained a braid $b\in \hat b$
such that
$\Lambda_{tr}(b)< \frac{2}{\pi}h(\hat{b})+\varepsilon$
for the a priori given small positive number $\varepsilon$. For each integer number $l$ the equalities $h(\widehat{b \Delta_3^{2l}})=h(\hat{b})$ and $\Lambda_{tr}(b \Delta_3^{2l})=\Lambda_{tr}(b)$ hold.
Hence,
the number of conjugacy classes $\hat{\mathbold{b}}$ of $\mathcal{B}_3\diagup\mathcal{Z}_3$ of positive entropy not exceeding $Y$
does not exceed the number of elements  $\mathbold{b}\in \mathcal{B}_3\diagup\mathcal{Z}_3 $ with $\Lambda_{tr}(\mathbold{b})< \frac{2}{\pi} Y +\varepsilon$. In other words,
$N^{entr}_{\mathcal{B}_3}(Y)
 \leq N^{\Lambda}_{\mathcal{B}_3}(\frac{2}{\pi} Y +\varepsilon)$. Since for $\varepsilon$ we may take any a priory given positive number, Theorem \ref{thm10.2} implies
\begin{equation}\label{eq10.32}
N^{entr}_{\mathcal{B}_3}(Y)
 \leq N^{\Lambda}_{\mathcal{B}_3}(\frac{2}{\pi} Y) \leq 4 e^{12 Y}\,.
\end{equation}
We obtained the upper bound. \hfill $\Box$

\medskip

For obtaining the lower bound we need the following preparations.

\begin{lemm}\label{lem10.20}
Suppose $b_1$ and $b_2$ are elements of the free group $\mathbb{F}_n$ in $n$ generators, that are both represented by cyclically reduced words. Then $b_1 $ and $b_2$ are conjugate if and only if
the word representing $b_2$ is obtained from the word representing $b_1$ by a cyclic permutation of terms.
\end{lemm}

\noindent {\bf Proof.}
It is enough to prove the following statement. If under the conditions of the lemma
$b_2=w^{-1} b_1 w$ for an element $ w\in \mathbb{F}_2$,
$w\neq\mbox{Id}$, then $b_2=w'^{-1}b'_1 w'$ where $b'_1$ is represented by a cyclically reduced word that is obtained from the reduced word representing $b_1$ by a cyclic permutation of terms, and $w'$ is represented by a reduced word which has less terms than the word representing $w$.

This statement is proved as follows. Write $w=w_1 w'$ where $w_1\in \mathbb{F}_n$ is represented by the first term of the word representing $w$, and $b_1= a' B_1' a''$ where $a'$ and $a''$ are represented by the first and the last term, respectively, of the word
that represents $b_1$. Then $b_2=w'^{-1}w_1^{-1} a' B_1' a'' w_1 w'$. If both relations $w_1^{-1} a'\neq {\rm Id}$ and $a'' w_1\neq {\rm Id}$ were true, then the first and the last term of the  reduced word representing $b_2$ would be a power of the same generator of $\mathbb{F}_n$, which contradicts the fact that $b_2$ can be represented by a cyclically reduced word. If either $w_1^{-1} a'= {\rm Id}$ or $a'' w_1= {\rm Id}$, then the reduced words representing $b_1'\stackrel{def}=w_1^{-1} a' B_1' a'' w_1$ and $b_1$ are cyclic permutations of each other. Hence, the statement is true.
\hfill $\Box$

\begin{lemm}\label{lem10.21} The following equalities hold.
\begin{align}\label{eq10.16a}
\Delta_3 \sigma_1 = & \sigma_2 \Delta_3\,, \nonumber\\
\Delta_3 \sigma_2 = & \sigma_1 \Delta_3\,, \nonumber\\
\sigma_1^{-1}(\sigma_2^{-4}\Delta_3^4)\sigma_1= & \sigma_2^2 \sigma_1^2 \sigma_2^2 \sigma_1^2\,,\nonumber\\
\sigma_2^{-1}(\sigma_1^{-4}\Delta_3^4)\sigma_2= & \sigma_1^2 \sigma_2^2 \sigma_1^2 \sigma_2^2\,.
\end{align}
\end{lemm}

\noindent{\bf Proof.}
The third equality is obtained as follows
\begin{align}\label{eq10.16b}
& \sigma_1^{-1} (\sigma_2^{-4}\Delta_3^4)\sigma_1
=\sigma_1^{-1} \sigma_2^{-1}\Delta_3^2\sigma_2^{-2}\Delta_3^2\sigma_2^{-1}\sigma_1\nonumber\\
=&(\sigma_1^{-1}\sigma_2^{-1})(\sigma_2 \sigma_1 \sigma_2)(\sigma_2 \sigma_1 \sigma_2)\sigma_2^{-1}\sigma_2^{-1}(\sigma_2 \sigma_1 \sigma_2)(\sigma_2 \sigma_1 \sigma_2)\sigma_2^{-1}\sigma_1\nonumber\\
=&\sigma_2^2 \sigma_1^2 \sigma_2^2 \sigma_1^2\,.
\end{align}
The fourth equality is obtained by conjugating the third equality by $\Delta_3$. \hfill $\Box$
\medskip

Recall that each braid $\beta\in \mathcal{B}_3$ can be written uniquely in the form $\beta=\sigma_j^k\beta_1 \Delta_3^{\ell}$ for an integer number $\ell$, a number $j=1,2,$ a number $k=0,1,$ and a braid $\beta_1$
that can be written as a reduced word in $\sigma_1^2$ and $\sigma_2^2$.
Consider two elements $\mathbold{b}_1$ and $\mathbold{b}_2$ of  $\mathcal{PB}_3\diagup \mathcal{Z}_3 \subset \mathcal{B}_3\diagup \mathcal{Z}_3$. We identify $\mathcal{PB}_3\diagup \mathcal{Z}_3$ with the free group in two generators $a_1$ and $a_2$.
\begin{lemm}\label{lem10.22} Suppose both, $\mathbold{b}_1$ and $\mathbold{b}_2$, are of the form
\begin{equation}\label{eq10.15c}
a_1^{\pm 2}\, a_2 ^{\pm 2}\cdots a_1^{\pm 2}\, a_2 ^{\pm 2}\;\;\;\;\mbox{or}\;\;\;\; a_2^{\pm 2}\, a_1 ^{\pm 2}\cdots a_2^{\pm 2}\, a_1 ^{\pm 2}\,,
\end{equation}
with a positive number of terms,
and the reduced word representing $\mathbold{b}_1$ has at least four terms.
Then $\mathbold{b}_2$ cannot be conjugated to $\mathbold{b}_1$ by an element $\mathbold{\beta}$ of  $\mathcal{B}_3\diagup \mathcal{Z}_3$ that can be represented by a braid
$\beta=\sigma_j\beta_1 \Delta_3^{\ell}$
with some $j$ and $\ell$ and $\beta_1$ being a word in $\sigma_1^2$ and $\sigma_2^2$.
\end{lemm}

\noindent{\bf Proof.}
Indeed, suppose the contrary,
\begin{equation}\label{eq10.15a}
\mathbold{b}_1=\mathbold{\beta}^{-1} \mathbold{b}_2 \mathbold{\beta}
\end{equation}
with $\mathbold{\beta}$ represented by $\beta=\sigma_j\beta_1 \Delta_3^{\ell}$,
where  $\ell=0,1,$ and $\beta_1\in \mathcal{PB}_3$ is a word in $\sigma_1^2$ and $\sigma_2^2$. For $j=1,2,$ we let $b_j$ be the representative of $\mathbold{b}_j$ which can be written as reduced word in $\sigma_1^2$ and $\sigma_2^2$. By equation \eqref{eq10.16a} there is an integer number $n$ such that the braid $b_2'\stackrel{def}=\sigma_j^{-1} b_2 \sigma_j \Delta_3^{2n}$ is a product of a positive even number of factors which either have alternately the form $\sigma_1^{\pm 4}$ and $(\sigma_2^2 \sigma_1^2 \sigma_2^2 \sigma_1^2)^{\pm 1}$, or they have alternately the form $\sigma_2^{\pm 4}$ and $(\sigma_1^2 \sigma_2^2 \sigma_1^2 \sigma_2^2)^{\pm 1}$. The braid ${\beta_1} \Delta_3^{\ell} b_1 \Delta_3^{-\ell}\beta_1^{-1}$ can also be written as reduced word in $\sigma_1^2$ and $\sigma_2^2$. Hence, by equation \eqref{eq10.15a} the two braids $b_2'$ and ${\beta_1} \Delta_3^{\ell} b_1 \Delta_3^{-\ell}\beta_1^{-1}$ must be equal.

Identify $\sigma_1^2$ with the generator $a_1$ and $\sigma_2^2$ with the generator $a_2$ of the free group $\mathbb{F}_2$.
The braid $b'_2$ is the product of at least two factors equal to $a_1^{\pm 2}$ and $(a_1a_2a_1a_2)^{\pm 1}$ alternately, or it is the product of at least two factors equal to $a_2^{\pm 2}$ and $(a_2a_1a_2a_1)^{\pm 1}$ alternately.
Hence, the reduced word representing $b'_2$ contains at least $4$ terms, and out of four consecutive terms at least two terms appear with power $+1$ or $-1$. 

Indeed, put $A=a_1a_2a_1a_2$. Replace in the product $A^{\ell_1}a_1^{\pm 2} A^{\ell_2}$ with $\ell_1=0,\pm 1,$ $\ell_2=0,\pm 1,$ each factor $A^{\ell_1}$ by the reduced word in $a_1$, $a_2$, representing it. We obtain a (possibly not reduced) word in $a_1$, $a_2$. 
If $\ell_1=0$ or $\ell_1=1$, and $\ell_2=0$ or $-1$, the word is already reduced.
If $\ell_1=-1$, and $\ell_2=0$ or $-1$,
the reduced word representing the product, consists of the first three letters of the word representing $A^{-1}=A^{\ell_1}$, a non-trivial power of $a_1$ and all letters of the word representing  $A^{\ell_2}$.
If $\ell_1=0$ or $\ell_1=1$, and $\ell_2=1$,
the reduced word representing the product, consists of all letters of the word representing $A^{\ell_1}$, a non-trivial power of $a_1$ and the last three letters of the word representing  $A= A^{\ell_2}$. Finally , if $\ell_1=-1$ and $\ell_2=1$ then the reduced word representing the product, consists of the first three letters of the word representing $A^{-1}=A^{\ell_1}$, a non-trivial power of $a_1$ and the last three letters of the word representing $A=A^{\ell_2}$.
Hence, the reduced word representing $b_2'$ contains for each factor $A^{\ell_j}$ in the whole product at least the two middle letters of the word representing the $A^{\ell_j}$. 
The case, when $b_2'$ is the product of $a_2a_1a_2a_1$ and $a_2^2$ is obtained by replacing the role of the generators $a_1$ and $a_2$.

Since the braid $\Delta_3^{\ell} b_1 \Delta_3^{-\ell}$  can be represented by a word of the form \eqref{eq10.15c} with at least four terms,
we arrived at the following conjugation problem in the free group $\mathbb{F}_2$. We have an element $B_1\in \mathbb{F}_2$ of the from \eqref{eq10.15c} with at least $4$ terms, and an element $B_2\in \mathbb{F}_2$ written as reduced word with at least $4$ terms, such that in each tuple of four consecutive terms
there are at least two terms that appear with power $+1$ or $-1$. Then there is no element of $\mathbb{F}_2$ that conjugates $B_1$ to $B_2$.

To prove the statement we assume the contrary.
Suppose the reduced word $a_{j_1}^{k_1}\ldots a_{j_{\ell}}^{k_{\ell}}$ representing $B_2$ is not cyclically reduced. Identify each term $ a_{j_{\ell'}}^{k_{\ell'}}$ of the word with the element of $\mathbb{F}_2$ represented by it.
The conjugate $(a_{j_1}^{k_1})^{-1}   B_2 a_{j_1}^{k_1}$ can be represented by a reduced word
for which the consecutive sequence of all terms except perhaps the last
is also a consecutive sequence of terms of the reduced word representing $B_2$. Continue consecutively in this way
with the $ a_{j_{\ell'}}^{k_{\ell'}}$ until we arrive at an element $B_2'\in \mathbb{F}_2$ that can be represented by a cyclically reduced word. The sequence of all consecutive terms except perhaps the last of this cyclically reduced word is also a sequence of consecutive terms
of the reduced word representing $B_2$.

Since $B_2'$ is conjugate to $B_1$ and both are represented by cyclically reduced words, the representing words have the same number of terms  (see Lemma \ref{lem10.20} ). By the assumption for $B_1$ the number of terms is at least $4$. Then the words $B_2'$ and $B_2$ have at least $3$ consecutive terms in common. Therefore $B_2'$ contains a term that appears with power $\pm 1$ which is a contradiction. The lemma is proved. \hfill $\Box$

\medskip

\noindent {\bf Proof of the lower bound of Theorem \ref{thm10.0}.}
Consider the elements of $\mathbold{b}\in\mathcal{PB}_3\diagup \mathcal{Z}_3 \subset \mathcal{B}_3\diagup \mathcal{Z}_3$ which can be represented by words of the form
\begin{equation}\label{eq10.15}
w=a_1^{\pm 2}\, a_2 ^{\pm 2}\cdots a_1^{\pm 2}\, a_2 ^{\pm 2}
\end{equation}
in the generators $a_1$ and $a_2$ of $\mathcal{PB}_3\diagup \mathcal{Z}_3$ with at least $4$ terms.
We denote the number of syllables (in this case the number of the terms $a_i^{k_i}$) of the word \eqref{eq10.15} by $2j$.
Since each word $w$ of form \eqref{eq10.15} is cyclically syllable reduced,  Theorem  \ref{thm10.1a} implies that the entropy $h(\hat{\mathbold{b}})$ of the conjugacy class of the element $\mathbold{b}\in \mathcal{PB}_3\diagup \mathcal{Z}_3$ represented by $w$ satisfies the inequality
\begin{align}\label{eq10.16}
h(\hat{\mathbold{b}}) \leq 300 \frac{\pi}{2} \mathcal{L}_+(w)=\frac{\pi}{2} \cdot 300 \cdot 2j \cdot \log 8\,.
\end{align}
Take $j=j'_0 \stackrel{def}= [\frac{Y}{300 \pi \log 8}]$.
For this choice of $j'_0$ the inequality  $\frac{\pi}{2}  \cdot 300 \mathcal{L}_+(w) \leq Y$ holds. Since we required that
the number of syllables $2j'_0$ is at least $4$, we get the condition $Y \geq 600 \pi \log 8$.
The number of different words of such kind equals $2^{2 j'_0}$.

We prove now that the number of different conjugacy classes of $\mathcal{B}_3\diagup\mathcal{Z}_3$ that can be represented by elements in $\mathcal{PB}_3$
corresponding to words \eqref{eq10.15} with $2j'_0$ syllables is not smaller than $\frac{2^{2j'_0}}{2j'_0}$. It is enough to prove the following claim. For each element
$\mathbold{b} \in \mathcal{PB}_3\diagup\mathcal{Z}_3$ of form \eqref{eq10.15} the number of
elements of $\mathcal{PB}_3\diagup\mathcal{Z}_3$ of form \eqref{eq10.15} that are conjugate to $\mathbold{b}$ by an element of $\mathcal{B}_3\diagup\mathcal{Z}_3$ does not exceed $2j'_0$.

Suppose two elements $\mathbold{b}_1$ and $\mathbold{b}_2$ of  $\mathcal{PB}_3\diagup \mathcal{Z}_3 \subset
\mathcal{B}_3\diagup \mathcal{Z}_3$ are represented by a word of form \eqref{eq10.15} and belong to the same conjugacy class, i.e. $\mathbold{b}_2=\mathbold{\beta} \mathbold{b}_1 \mathbold{\beta}^{-1}$ for an element $\mathbold {\beta} \in \mathcal{B}_3\diagup \mathcal{Z}_3 $. Then by Lemma \ref{lem10.22} the element $\mathbold{\beta}$ is represented by an element $\beta= \beta_1 \Delta_3^{\ell}$ in $\mathcal{B}_3$ with $\beta_1$ being a word in $\sigma_1^2$ and $\sigma_2^2$ and $\ell=0,1$. Put $\mathbold{b}'_1= {\mathbold{\Delta}_3}^{\ell} \mathbold{b}_1{\mathbold{\Delta}_3}^{-\ell}$ for the element ${\mathbold{\Delta}_3}=\Delta_3\diagup \mathcal{Z}_3$. Then $\mathbold{b}'_1$ is represented by a word of form \eqref{eq10.15c}.
If $\ell=0$ then the first terms of the reduced words representing $\mathbold{b}'_1$ and $\mathbold{b}_1$ are powers of the same generator $a_j$, if $\ell=1$ then the first terms of  the reduced words representing $\mathbold{b}'_1$ and $\mathbold{b}_1$ are powers of different generators. By
Lemmas \ref{lem10.20} and \ref{lem10.21} the reduced word representing the element $\mathbold{\beta}_1\mathbold{b}'_1 \mathbold{\beta}_1^{-1}$ is obtained from  the reduced word representing $\mathbold{b}'_1$ by a cyclic permutation of terms.
Thus, either the reduced word representing $\mathbold{b}_2$ is obtained by an even cyclic permutation of terms of the reduced  word representing $\mathbold{b}_1$, or by an odd cyclic permutation of the terms of the reduced word representing ${\mathbold{\Delta}_3} \mathbold{b}_1 {\mathbold{\Delta}_3}^{-1}$. Hence the number of different elements of $\mathcal{PB}_3\diagup \mathcal{Z}_3$ that can be represented by words of form  \eqref{eq10.15}
and are obtained from $\mathbold{b}_1$ by conjugation with an element of $\mathcal{B}_3\diagup \mathcal{Z}_3$ does not exceed the number of cyclic permutations of $2j'_0$ letters, i.e. $2j'_0$. (Notice, that e.g. by cyclically permuting the terms of a word with symmetries we may sometimes arrive at the same word.)

We proved that the number of different conjugacy classes of $\mathcal{PB}_3\diagup \mathcal{Z}_3$ represented by words of form  \eqref{eq10.15} is not smaller than $\frac{2^{2j'_0}}{2j'_0} $.

We obtain
\begin{equation}\label{eq10.17}
N^{entr}_{\mathcal{B}_3}(Y) \geq \frac{2^{2j'_0}}{2j'_0}\,.
\end{equation}
Notice that $\frac{2^j}{j}\geq 2 $ for natural $j$. Indeed, the function $x \to \frac{2^x}{x}$ increases for $x \geq 2$ (since $ (\frac{2^x}{x})' = -\frac{2^x}{x^2}+\frac{2^x \log 2}{x}$ , and $\log 2 >0.6$) and for $j=1$ and $2$ the expression equals $2$).
Hence, $\frac{2^{2j}}{2j}\geq 2^{j}$. Hence, for $Y \geq 600 \pi \log 8$
\begin{equation}\label{eq10.18}
N^{entr}_{\mathcal{B}_3}(Y) \geq  2^{j'_0}\geq \frac{1}{2} 2^{\frac{Y}{300 \pi \log 8}}
= \frac{1}{2} \exp(\frac{Y \log 2 }{300 \pi \log 8})
= \frac{1}{2} \exp(\frac{Y}{900 \pi}).
\end{equation}
The lower bound of the theorem is proved. \hfill $\Box$

\medskip
The author is grateful to the Max-Planck-Institute where the work on the paper was started, and to the IHES where the paper was finished. She would also like to thank the referee for the careful reading and for pointing out a mistake in the end of the first version.

\end{document}